\documentclass[11pt]{article}
\usepackage{stmaryrd}
\usepackage{mathrsfs}
\usepackage[centertags]{amsmath}
\usepackage{amsfonts, dsfont}
\usepackage{amssymb}
\usepackage{amsthm}

\input xy
\xyoption{all}

\theoremstyle{plain}
\newtheorem{thm}{Theorem}[section]
\newtheorem{cor}[thm]{Corollary}
\newtheorem{lem}[thm]{Lemma}
\newtheorem{prop}[thm]{Proposition}

\theoremstyle{definition}
\newtheorem{defn}[thm]{Definition}

\newtheorem*{remarks}{\textbf{Remarks}}
\newtheorem*{ack}{Acknowledgments}

\newcommand{\bd}{\begin{defn}}
\newcommand{\ed}{\end{defn}}
\newcommand{\bl}{\begin{lem}}
\newcommand{\el}{\end{lem}}
\newcommand{\bp}{\begin{prop}}
\newcommand{\ep}{\end{prop}}
\newcommand{\bt}{\begin{thm}}
\newcommand{\et}{\end{thm}}
\newcommand{\bc}{\begin{cor}}
\newcommand{\ec}{\end{cor}}
\newcommand{\br}{\begin{remarks}}
\newcommand{\er}{\end{remarks}}
\newcommand{\ba}{\begin{array}}
\newcommand{\ea}{\end{array}}
\newcommand{\bpf}{\begin{proof}}
\newcommand{\epf}{\end{proof}}

\newcommand{\Z}{\mathds{Z}}
\newcommand{\Q}{\mathds{Q}}
\newcommand{\Zp}{\mathds{Z}_{p}}
\newcommand{\Qp}{\mathds{Q}_{p}}

\newcommand{\Ga}{\Gamma}

\newcommand{\La}{\Lambda}

\newcommand{\m}{\mathfrak{m}}

\DeclareMathOperator{\Gal}{Gal} \DeclareMathOperator{\Hom}{Hom}
\DeclareMathOperator{\Ext}{Ext} 
 
\newcommand{\ot}{\otimes}

\newcommand{\ilim}{\displaystyle \mathop{\varinjlim}\limits}
\newcommand{\plim}{\displaystyle \mathop{\varprojlim}\limits}

\newcommand{\cts}{\mathrm{cts}}

\newcommand{\cyc}{\mathrm{cyc}}

\newcommand{\lra}{\longrightarrow}

\newcommand{\sbs}{\subseteq}

\newcommand{\ps}[1]{\llbracket #1 \rrbracket}

\begin{document}

\title{On the pseudo-nullity of the dual fine Selmer groups}
\author{Meng Fai Lim}
\date{}
\maketitle

\begin{abstract} \footnotesize
\noindent  In this paper, we will study the pseudo-nullity of the
dual fine Selmer group and its related question. We investigate
certain situations, where one can deduce the pseudo-nullity of the
dual fine Selmer group of a general Galois module over a admissible
$p$-adic Lie extension $F_{\infty}$ from the knowledge of the
pseudo-nullity of the Galois group of the maximal abelian unramified
pro-$p$ extension of $F_{\infty}$ at which every prime of
$F_{\infty}$ above $p$ splits completely. In particular, this gives
us a way to construct examples of the pseudo-nullity of the dual
fine Selmer group of a Galois module that is unramified outside $p$.
We will illustrate our results with many examples.

\medskip
\noindent Keywords and Phrases: Fine Selmer groups, $S$-admissible
$p$-adic Lie extensions, pseudo-null.

\smallskip \noindent
2010 Mathematics Subject Classification: Primary 11R23; Secondary
11R34, 11F80, 11S25.

\end{abstract}

\normalsize
\section{Introduction} \label{introduction}

Throughout the paper, $p$ will denote a fixed odd prime. Let $F$ be
a number field. Let $S$ be a finite set of primes of $F$ which
contains the primes above $p$ and the infinite primes. We then
denote by $F_S$ the maximal algebraic extension of $F$ which is
unramified outside $S$. For any algebraic (possibly infinite)
extension $\mathcal{L}$ of $F$ contained in $F_S$, we write
$G_S(\mathcal{L}) = \Gal(F_S/\mathcal{L})$. Let $R$ be a commutative
complete regular local ring with maximal ideal $\m$ and residue
field $k$, where $k$ is finite of characteristic $p$. We denote by
$T$ a finitely generated free $R$-module with a continuous
$R$-linear $G_S(F)$-action. Here $T$ is endowed with the canonical
topology arising from its filtration by the powers of $\m$.

Let $v$ be a prime in $S$. For every finite extension $L$ of $F$
contained in $F_S$, we define
 \[ K_v^2(T/L) = \bigoplus_{w|v}H^2(L_w, T),\]
where $w$ runs over the (finite) set of primes of $L$ above $v$. If
$\mathcal{L}$ is an infinite extension of $F$ contained in $F_S$, we
define
\[ \mathcal{K}_v^2(T/\mathcal{L}) = \plim_L K_v^2(T/L),\]
where the inverse limit is taken over all finite extensions $L$ of
$F$ contained in $\mathcal{L}$. For any algebraic (possibly
infinite) extension $\mathcal{L}$ of $F$ contained in $F_S$, the
\textit{dual fine Selmer group} of $T$ over $\mathcal{L}$ (with
respect to $S$) is defined to be
\[ Y_S(T/\mathcal{L}) = \ker\Big(H^2_S(\mathcal{L}/F, T)\lra \bigoplus_{v\in S} \mathcal{K}_v^2(T/\mathcal{L})
\Big), \] where we write $H^2_S(\mathcal{L}/F, T) =
\plim_{L}H^2(G_S(L), T)$.

Before continuing, we first compare our definition with that in
\cite{CS, Jh, JhS, Lim}. Let $W = \Hom_{\cts}(T, \mu_{p^{\infty}})$,
where $\mu_{p^{\infty}}$ denotes the group of all $p$-power roots of
unity. Then as loc. cit., the fine Selmer group of $W$ over
$\mathcal{L}$ (with respect to $S$) is defined to be
\[ R_S(W/\mathcal{L}) =
\ker\Big(H^1(G_S(\mathcal{L}), W)\lra \ilim_L\bigoplus_{w|v}H^1(L_w,
W)\Big).
\] Here the direct limit is taken over all finite extensions $L$ of
$F$ contained in $\mathcal{L}$. Then by the Poitou-Tate duality, we
have that $Y_S(T/\mathcal{L})$ is the Pontryagin dual of
$R_S(W/\mathcal{L})$ (and hence the terminology ``dual fine Selmer
group").

To facilitate further discussion, we recall the following conjecture
which was studied in \cite{CS, JhS, Lim}.

\medskip \noindent \textbf{Conjecture A :} For
any number field $F$,\, $Y_S(T/F^{\cyc})$ is a finitely generated
$R$-module, where $F^{\cyc}$ is the cyclotomic $\Zp$-extension of
$F$.

\medskip
We say that $F_{\infty}$ is an \textit{$S$-admissible $p$-adic Lie
extension} of $F$ if (i) $\Gal(F_{\infty}/F)$ is a compact pro-$p$
$p$-adic Lie group, (ii) $F_{\infty}$ contains $F^{\cyc}$ and (iii)
$F_{\infty}$ is contained in $F_S$. 
Write $G=\Gal(F_{\infty}/F)$ and $H=\Gal(F_{\infty}/F^{\cyc})$. We
recall that a finitely generated $R\ps{G}$-module $M$ is said to be
\textit{torsion} (resp., \textit{pseudo-null}) if $
\Ext_{R\ps{G}}^{i}(M, R\ps{G}) = 0$ for $i=0$ (resp., $i=0,1$). If
the $R\ps{G}$-module $M$ in question is finitely generated over
$R\ps{H}$, then it follows from a well-known result of Venjakob
(cf.\ \cite{V}) that it is a pseudo-null $R\ps{G}$-module if and
only if it is a torsion $R\ps{H}$-module.

Now if Conjecture A is valid, it follows from a similar argument to
that in \cite[Lemma 3.2]{CS} that $Y_S(T/F_{\infty})$ is finitely
generated over $R\ps{H}$. We can then state the following question,
noting the above result of Venjakob.

\medskip

\noindent \textbf{Question B :}  Let $F_{\infty}$ be an
$S$-admissible $p$-adic Lie extension of $F$ of dimension $>1$, and
suppose that $Y_S(T/F^{\cyc})$ is a finitely generated $R$-module.
Is $Y_S(T/F_{\infty})$ a pseudo-null $R\ps{G}$-module, or
equivalently a torsion $R\ps{H}$-module?

\medskip  When $T$ is
the Tate module of all the $p$-power roots of unity, the dual fine
Selmer group is precisely $\Gal(K(F_{\infty})/F_{\infty})$, where
$K(F_{\infty})$ is the maximal abelian unramified pro-$p$ extension
of $F_{\infty}$ at which every prime of $F_{\infty}$ above $p$
splits completely. In this context, Hachimori and Sharifi \cite{HS}
has constructed a class of admissible $p$-adic Lie extensions
$F_{\infty}$ of $F$ of dimension $>1$ such that
$\Gal(K(F_{\infty})/F_{\infty})$ is not pseudo-null. Despite so,
they have speculated that the pseudo-nullity statement should hold
for admissible $p$-adic extensions that ``come from algebraic
geometry'' (see \cite[Question 1.3]{HS} for details, and see also
\cite[Conjecture 7.6]{Sh2} for a related assertion and \cite{Sh3}
for positive results in this direction).

When $T$ is the Tate module of an elliptic curve, this is precisely
\cite[Conjecture B]{CS}.  In this context, the conjecture has also
been studied and verified in certain numerical examples (see
\cite{Bh, Jh, Lim, Oc}). When $T$ is the $R(1)$-dual of the Galois
representation attached to a normalized eigenform ordinary at $p$,
this is \cite[Conjecture B]{Jh}. In the case when $T$ is the
$R(1)$-dual of the Galois representation coming from a $\La$-adic
form, this is \cite[Conjecture 2]{Jh}. We should mention that in the
formulation of their conjectures in \cite{CS, Jh}, they did not have
any restriction on the admissible $p$-adic Lie extensions.

The following question is the theme of this paper.

\medskip

\noindent \textbf{Question B$'$:}  Let $F_{\infty}$ be an
$S$-admissible $p$-adic Lie extension of a number field $F$ of
dimension $>1$ with the property that $G_S(F_{\infty})$ acts
trivially on $T/\m T$. Suppose that $\Gal(K(F_{\infty})/F_{\infty})$
is a finitely generated $R\ps{H}$-module. Can one deduce that
$Y_S(T/F_{\infty})$ is a pseudo-null $R\ps{G}$-module from the
knowledge that $\Gal(K(F_{\infty})/F_{\infty})$ is a pseudo-null
$\Zp\ps{G}$-module?

\medskip We note that the assumptions in the above question imply that
$Y_S(T/F_{\infty})$ is a finitely generated $R\ps{H}$-module (see
\cite[Theorem 3.5, Lemma 5.2]{Lim}). In fact, one can think of
Question B$'$ as whether the pseudo-null analogue of \cite[Theorem
3.5]{Lim} is valid. In this paper, we will be studying Question
B$'$. Namely, we will give some results which partially answer
Question B$'$ (see Theorems \ref{pseudo-null main} and
\ref{pseudo-null main2}). We then apply Theorem \ref{pseudo-null
main} to establish the pseudo-nullity of the dual fine Selmer group
of a certain class of Galois modules that are unramified outside
$p$, and therefore, this gives a positive answer to the question B
for this particular class of Galois modules.
We then give several examples of these pseudo-null dual fine Selmer
groups which come from elliptic curves, modular forms and abelian
varieties.

We end the introduction by mentioning that in view of the
counterexamples of Hachimori-Sharifi \cite{HS} and the
pseudo-nullity conjectures made in \cite{CS, Jh}, it seems unlikely
that the converse of the statement of Question B$'$ will hold. It
will be interesting to find a numerical counterexample to the
converse statement which we are not able to at the moment.

\section{Question B$'$}

We retain the notion and notation of Section \ref{introduction}. We
now record a lemma which gives a relationship between the dual fine
Selmer group of $T$ over $F_{\infty}$ and the Galois group
$\Gal(K(F_{\infty})/F_{\infty})$ under the assumption that
$G_S(F_{\infty})$ acts trivially on $T$. We denote $M(n)$ to be the
$n$th Tate twist of $M$.

\bl \label{fine selmer and class group} Suppose that $F$ contains
$\mu_p$. Let $F_{\infty}$ be an $S$-admissible $p$-adic Lie
extension of $F$ such that $G_S(F_{\infty})$ acts trivially on $T$.
Then we have an isomorphism
\[ Y_S(T/F_{\infty}) \cong \Gal(K(F_{\infty})/F_{\infty}) \ot_{\Zp} T(-1) .\]
 \el

\bpf For a pro-$p$ group or a discrete $p$-primary group $M$, its
Pontryagin dual is denoted to be $M^{\vee}= \Hom_{\cts}(M,
\Qp/\Zp)$. As the group $G_S(F_{\infty})$ acts trivially on $T$ and
$\mu_{p^{\infty}} \sbs F_{\infty}$, it also acts trivially on
$T^{\vee}(1)$. Therefore, one calculates that
\[ Y_S(T/F_{\infty})^{\vee} = \Hom_{\Zp}\big(\!\Gal(K(F_{\infty})/F_{\infty}), T^{\vee}(1)\big). \]
 On the other hand, one has the following adjunction
isomorphism
\[ \Big(\Gal(K(F_{\infty})/F_{\infty}) \ot_{\Zp}T(-1)\Big)^{\vee}
\cong \Hom_{\Zp}\big(\!\Gal(K(F_{\infty})/F_{\infty}),
T^{\vee}(1)\big).\] Combining this with the above equality and
taking Pontryagin dual, we obtain the required isomorphism. \epf

As we are mostly interested in $S$-admissible $p$-adic Lie
extensions $F_{\infty}$ that satisfy the following condition, we
shall give a name to it.

\medskip
\noindent \textbf{(Dim$_S$) :}  For each $v\in S$, the decomposition
group of $G = \Gal(F_{\infty}/F)$ at $v$ has dimension $\geq 2$.

\medskip
We record another result (cf.\ \cite[Theorem 5.4]{Lim}, see also
\cite[Theorem 10]{Jh}) which is the main ingredient in proving our
main results. In view of the discussion in this paper, we will state
a slightly strengthened version of \cite[Theorem 5.4]{Lim}.

\bp \label{pseudo-null tech}
 Let $F_{\infty}$ be an $S$-admissible
$p$-adic Lie extension of $F$ which satisfies $(\mathbf{Dim}_S)$.
Assume that $Y_S(T/F_{\infty})$ is a finitely generated
$R\ps{H}$-module.  Suppose that there exists a prime ideal
$\mathfrak{p}$ of $R$ such that the ring $R/\mathfrak{p}$ is also
regular local. If $Y_S\big((T/\mathfrak{p}T)/F_{\infty}\big)$ is a
pseudo-null $R/\mathfrak{p}\ps{G}$-module, then $Y_S(T/F_{\infty})$
is a pseudo-null $R\ps{G}$-module. \ep

\bpf
 Let $F_0$ be a finite
extension of $F$ contained in $F_{\infty}$ such that
$G_0:=\Gal(F_{\infty}/F_0)$ is a compact pro-$p$ $p$-adic Lie group
without $p$-torsion. Clearly, the decomposition group of $G_0$ at
$v$ for every $v\in S$ has dimension $\geq 2$. Write
$H_0=\Gal(F_{\infty}/F_0^{\cyc})$. Since $F_0$ is a finite
$p$-extension of $F$, $H_0$ is a subgroup of $H$ with finite index.
Therefore, we also have that $Y_S(T/F_{\infty})$ is a finitely
generated $R\ps{H_0}$-module. Now it is easy to verify that (or see
\cite[Proposition 5.4.17]{NSW}), for every $i\geq 0$, one has an
isomorphism
\[ \Ext_{R\ps{G}}^{i}(M, R\ps{G}) \cong
\Ext_{R\ps{G_0}}^{i}(M, R\ps{G_0})\] for any finitely generated
$R\ps{G}$-module $M$. Hence, we are reduced to proving the
proposition over $G_0$. But this is precisely the statement of
\cite[Theorem 5.4]{Lim}. \epf

 We can now proceed to prove the following two theorems which
 partially answer Question B$'$.

 \bt \label{pseudo-null main}
 Suppose that $F$ contains $\mu_p$. Let $F_{\infty}$ be an $S$-admissible
$p$-adic Lie extension of $F$ which satisfies $(\mathbf{Dim}_S)$.
Assume that $G_S(F_{\infty})$ acts trivially on $T/\m T$, and assume
that $\Gal(K(F_{\infty})/F_{\infty})$ is a finitely generated
$\Zp\ps{H}$-module. If $\Gal(K(F_{\infty})/F_{\infty})/p$ is a
pseudo-null $\Z/p\Z\ps{G}$-module, then $Y_S(T/F_{\infty})$ is a
pseudo-null $R\ps{G}$-module. \et

Note that it follows from Proposition \ref{pseudo-null tech} that if
$\Gal(K(F_{\infty})/F_{\infty})/p$ is a pseudo-null
$\Z/p\Z\ps{G}$-module, then $\Gal(K(F_{\infty})/F_{\infty})$ is a
pseudo-null $\Zp\ps{G}$-module.

\bpf[Proof of Theorem \ref{pseudo-null main}] Firstly, as noted in
the introduction, since $G_S(F_{\infty})$ acts trivially on $T/\m T$
and $\Gal(K(F_{\infty})/F_{\infty})$ is a finitely generated
$\Zp\ps{H}$-module, it follows from \cite[Theorem 3.5]{Lim} that
$Y_S(T/F_{\infty})$ is a finitely generated $R\ps{H}$-module. Now by
Lemma \ref{fine selmer and class group}, we have an isomorphism
  \[ Y_S((T/\m T)/F_{\infty}) \cong \Gal(K(F_{\infty})/F_{\infty}) \ot_{\Zp} T/\m T(-1)  .\]
  Since $p$ kills $T/\m T$, the latter is equal to
  \[ (\Gal(K(F_{\infty})/F_{\infty})/p) \ot_{\Zp} T/\m T(-1). \]
It then follows from the pseudo-nullity hypothesis of the theorem
that $Y_S((T/\m T)/F_{\infty})$ is a pseudo-null
$R/\m\ps{\Gal(F_{\infty}/\Q(\mu_p))}$-module. Therefore, we may now
apply Proposition \ref{pseudo-null tech} to obtain the required
conclusion.
 \epf

\medskip
Before stating the next result, we introduce another notation. For a
prime ideal $\mathfrak{p}$ of $R$, we denote the residual
representation of $\rho$ mod $\mathfrak{p}$ to be
 \[ \rho_{\mathfrak{p}} : G_S(F) \lra
 \mathrm{Aut}_{R/\mathfrak{p}}(T/\mathfrak{p}T). \]

\bt \label{pseudo-null main2}
 Suppose that $F$ contains $\mu_p$. Let $F_{\infty}$ be an $S$-admissible
$p$-adic Lie extension of $F$ which satisfies $(\mathbf{Dim}_S)$.
Assume that $\Gal(K(F_{\infty})/F_{\infty})$ is a finitely generated
$\Zp\ps{H}$-module.  Suppose that there exists a prime ideal
$\mathfrak{p}$ of $R$ such that the ring $R/\mathfrak{p}$ is a
regular local ring which is a finite free $\Zp$-algebra and such
that $\Gal(F_S/F_{\infty})\sbs \ker \rho_{\mathfrak{p}}$.  If
$\Gal(K(F_{\infty})/F_{\infty})$ is a pseudo-null
$\Zp\ps{G}$-module, then $Y_S(T/F_{\infty})$ is a pseudo-null
$R\ps{G}$-module. \et

\bpf By Proposition \ref{pseudo-null tech}, one is reduced to
showing that $Y_S((T/\mathfrak{p}T)/F_{\infty})$ is a pseudo-null
$R/\mathfrak{p}\ps{G}$-module. Since $F$ contains $\mu_p$, the field
$F_{\infty}$ necessarily contains $\mu_{p^{\infty}}$, and therefore,
the group $G_S(F_{\infty})$ will act trivially on
$(T/\mathfrak{p}T)^{\vee}(1)$. By Lemma \ref{fine selmer and class
group}, we have
\[ Y_S((T/\mathfrak{p}T)/F_{\infty})\cong \Gal(K(F_{\infty})/F_{\infty})\ot_{\Zp}T/\mathfrak{p}T(-1). \]
Since the ring $R/\mathfrak{p}$ is a finite free $\Zp$-algebra and
$\Gal(K(F_{\infty})/F_{\infty})$ is a pseudo-null
$\Zp\ps{G}$-module, it follows that
$Y_S((T/\mathfrak{p}T)/F_{\infty})$ is a pseudo-null
$R/\mathfrak{p}\ps{G}$-module, as required. \epf

\section{Question B} \label{proof section}

In this section, we will apply Theorem \ref{pseudo-null main} to
prove the pseudo-nullity of the dual fine Selmer group of certain
classes of Galois modules that are unramified outside $p$ over an
$p$-adic admissible $p$-adic extension that is also unramified
outside $p$.

\textit{Therefore, from now on, we will say an admissible $p$-adic
Lie extension to mean an $S$-admissible $p$-adic Lie extension,
where the set $S$ consists \textbf{only} of the primes above $p$. We
will also assume that $T$ is unramified outside $p$.}

\bt \label{main}

Let $p$ be a regular prime. Suppose that $F_{\infty}$ is an
admissible $p$-adic Lie extension of $\Q(\mu_p)$ of dimension $>1$.
Consider the following statements.
            \begin{enumerate}
          \item[$(1)$] $G_S(F_{\infty})$ acts trivially on $T$,
          \item[$(2)$] $G_S(F_{\infty})$ acts trivially on $T/\m T$,
          \item[$(3)$] For the $($unique$)$ prime $v$ of $\Q(\mu_p)$ above $p$, the
          decomposition group of $\Gal(F_{\infty}/\Q(\mu_p))$ at $v$ has
          dimension $\geq 2$.
        \end{enumerate}
If statement $(1)$ holds, then $Y_S(T/F_{\infty}) =0$. If statements
$(2)$ and $(3)$ hold, then $Y_S(T/F_{\infty})$ is a pseudo-null
$R\ps{G}$-module. \et

The above theorem can be seen as a generalization of the results in
\cite{Oc}. The conclusion of the theorem under the assumption of
statement (1) is probably well-known, but nevertheless, we have
included it for completeness. In the case when $p$ is irregular, we
also have a statement in the spirit of the previous theorem, but
under more restrictive conditions.

\bt \label{main2}

Let $p$ be an irregular prime $<1000$. Suppose that $F_{\infty}$ is
an admissible $p$-adic Lie extension of $\Q(\mu_p)$ $($of dimension
$>1)$ which satisfies all of the following conditions.
 \begin{enumerate}
                    \item[$(1)$] $F_{\infty}$ is unramified outside $p$,
                \item[$(2)$] $G_S(F_{\infty})$ acts trivially on $T/\m
                T$,
                \item[$(3)$] $F_{\infty}$ contains $\Q(\mu_{p^{\infty}},
                p^{-p^{\infty}})$.
             \end{enumerate}
  Then $Y_S(T/F_{\infty})$ is a pseudo-null $R\ps{G}$-module. \et

In preparation for the proofs of the above two results, we collect
some preliminary results on the group
$\Gal(K(F_{\infty})/F_{\infty})$, where $K(F_{\infty})$ is the
maximal abelian unramified pro-$p$ extension of $F_{\infty}$ at
which every prime of $F_{\infty}$ above $p$ splits completely. We
begin stating the following well-known statement (for instance, see
\cite[Section 4]{Oc}).

\bp \label{main ingredient}
  Let $p$ be a regular prime. Let
$F_{\infty}$ be an admissible $p$-adic Lie extension of $\Q(\mu_p)$.
Then one has $\Gal(K(F_{\infty})/F_{\infty}) = 0$. \ep

The next result is concerned about the structure of group
$\Gal(K(F_{\infty})/F_{\infty})$ for irregular primes $<1000$ which
is an immediate consequence of \cite[Corollary 5.9]{Sh2} and
\cite[Propositions 3.3 and 2.1a]{Sh3}

\bp \label{main ingredient2}

Let $p$ be an irregular prime $<1000$. Let $F_{\infty}$ be an
admissible $p$-adic extension of $\Q(\mu_p)$ which contains
$\Q(\mu_{p^{\infty}}, p^{-p^{\infty}})$. Then
$\Gal(K(F_{\infty})/F_{\infty})$ is a finitely generated
$\Zp\ps{\Gal(F_{\infty}/\Q(\mu_{p^{\infty}},
p^{-p^{\infty}}))}$-module.  \ep

We can now prove our theorems.

\bpf[Proof of Theorem \ref{main}] Suppose that statement (1) holds.
Then by Lemma \ref{fine selmer and class group}, we have an
isomorphism
  \[ Y_S(T/F_{\infty}) \cong \Gal(K(F_{\infty})/F_{\infty}) \ot_{\Zp} T(-1)  .\]
By Proposition \ref{main ingredient}, this in turn implies that
$Y_S(T/F_{\infty}) =0$.

Now suppose that statements (2) and (3) hold. Then clearly
$\Gal(K(F_{\infty})/F_{\infty})/p$ is a pseudo-null
$\Z/p\Z\ps{G}$-module by Proposition \ref{main ingredient}. The
required conclusion is now immediate from an application of Theorem
\ref{pseudo-null main}.
 \epf

\bpf[Proof of Theorem \ref{main2}] It follows from Proposition
\ref{main ingredient2} that the group
$\Gal(K(F_{\infty})/F_{\infty})/p$ is a finitely generated
$\Z/p\Z\ps{\Gal(F_{\infty}/\Q(\mu_{p^{\infty}},
p^{-p^{\infty}}))}$-module and, in particular, a pseudo-null
$\Z/p\Z\ps{\Gal(F_{\infty}/\Q(\mu_p))}$-module. On the other hand,
it follows from \cite[Lemma 3.9]{HV} that \linebreak
$\Q(\mu_{p^{\infty}}, p^{-p^{\infty}})$, and hence $F_{\infty}$,
satisfies $(\mathbf{Dim}_S)$ for $S= \{p, \infty\}$. The required
conclusion now follows from Theorem \ref{pseudo-null main}.
 \epf

\section{Examples} \label{examples section}

In this section, we discuss some numerical examples of Theorems
\ref{main} and \ref{main2}.

(a) (See also \cite{Oc}) Let $E$ be an elliptic curve defined over
$\Q$ that has good reduction away from $p$ and possesses a
$\Q$-rational isogeny. A complete list (up to isogeny) of such
elliptic curves can be found in \cite[Table 2]{RT}. By
\cite[Proposition 5]{RT}, we have that $\Q(E[p])$ is a finite
$p$-extension of $\Q(\mu_p)$.  Now let $T$ be the Tate module of the
$p$-division points of such an elliptic curve. We first consider the
case that $p$ is regular (i.e., $p= 3,7, 11, 19, 43, 163$). It then
follows from Theorem \ref{main} that $Y_S(T/F_{\infty})$ is a
pseudo-null (resp., trivial)
$\Zp\ps{\Gal(F_{\infty}/\Q(\mu_p))}$-module when $F_{\infty}$ is an
admissible $p$-adic Lie extension containing $\Q(E[p])$ (resp.,
$\Q(E[p^{\infty}])$) whose decomposition group of
$\Gal(F_{\infty}/\Q(\mu_p))$ at the unique prime $v$ of $F$ has
dimension $\geq 2$.

Now consider the case when $p = 67$ (an irregular prime). By an
application of Theorem \ref{main2}, $Y_S(T/F_{\infty})$ is a
pseudo-null $\Zp\ps{\Gal(F_{\infty}/\Q(\mu_{67}))}$-module whenever
$F_{\infty}$ is an admissible $67$-adic Lie extension that contains
$\Q(E[67], \mu_{67^{\infty}}, 67^{-67^{\infty}})$. Examples of such
admissible $67$-adic Lie extensions are
\[\Q(\mu_{67^{\infty}}, 67^{-67^{\infty}}, E[67]), \quad\Q(E[67^{\infty}], 67^{-67^{\infty}}). \]

\medskip
(b) The discussion in (a) can also be applied to the case of a
modular form of higher weight with appropriate assumptions. Let $N$
be a power of the prime $p$ (we also allow $N=1$). Let $f \in
S_k(N)$ be a primitive cuspidal modular form of positive weight $k
\geq 2$, level $N$ and trivial character for the group $\Ga_0(N)$.
Fix an embedding $\Q$ in $\Qp$. By the results of Eichler, Shimura
and Deligne, there is an associated Galois representation, which we
denote by
 \[ \rho_f : \Gal(\bar{\Q}/\Q)\lra \mathrm{GL}_2(\Q_p) \]
which is unramified outside $p$ (since $N$ is a power of $p$). After
conjugation, one may assume that $\rho_f$ takes values in
$\mathrm{GL}_2(\Zp)$. Let $T$ be the Galois module associated to
$\rho_f$ which is a free $\Zp$-module of rank 2. Denote the residual
representation to be
\[ \bar{\rho}_f : \Gal(\bar{\Q}/\Q)\lra \mathrm{GL}_2(\Z/p\Z). \]
Let $K$ be the finite extension of $\Q$ corresponding to fixed field
of $\ker \bar{\rho}_f$. If the reduced representation $\bar{\rho}_f$
is reducible and contains $\mu_p$ as a sub-representation, then the
field $\Q(\mu_p)$ is contained in $K$ and the Galois group
$\Gal(K/\Q(\mu_p))$ is a finite $p$-group. Therefore, one can obtain
examples of psuedo-nullity of the dual fine Selmer group of such $T$
as in (a).

\medskip
(c) Assume that $p \geq 5$. Let $n\geq 1$. Let $a, b ,c$ be integers
such that $1\leq a, b,c < p^n$, $a+b+c = p^n$ and at least one of
(and hence at least two of) $a, b, c$ is not divisible by $p$. Let
$J=J_{a,b,c}$ be the Jacobian variety of the curve $y^{p^n} =
x^a(1-x)^b$. Set $T$ to be the Tate module of the $p$-division
points of $J$. Write $R = \Zp[\zeta_{p^n}]$ and $\pi =
1-\zeta_{p^n}$. Note that $\pi$ is a generator of the maximal ideal
of $R$. It follows from \cite[P.77]{Ih} and \cite[Chap.\ II, Theorem
5B]{Ih} (see also \cite[Section 2]{Gr} for the case $n=1$) that $T$
is a free $R$-module of rank one and that $G_S(\Q(\mu_{p^n}))$ acts
trivially on $T/\pi T$. Now if the prime $p$ is regular, it follows
from Theorem \ref{main} that $Y_S(T/F_{\infty})$ is a pseudo-null
$\mathcal{O}\ps{\Gal(F_{\infty}/\Q(\mu_p))}$-module, where
$F_{\infty}$ is an admissible $p$-adic Lie extension of $\Q(\mu_p)$
whose decomposition group of $\Gal(F_{\infty}/\Q(\mu_p))$ at the
unique prime $v$ of $F$ has dimension $\geq 2$. Now if $p$ is an
irregular prime $<1000$, one applies Theorem \ref{main2} to conclude
that $Y_S(T/F_{\infty})$ is a pseudo-null
$\mathcal{O}\ps{\Gal(F_{\infty}/\Q(\mu_p))}$-module for every
admissible $p$-adic Lie extension $F_{\infty}$ of $\Q(\mu_p)$ which
contains $\Q(\mu_{p^{\infty}}, p^{-p^{\infty}})$.

\begin{ack}
     This work was written up when the author was doing his postdoctoral at the GANITA Lab
    at the University of Toronto. He would like to acknowledge the
    hospitality and conducive working conditions provided by the GANITA
    Lab and the University of Toronto. The author would also
    like to thank Romyar Sharifi for his interest and
    comments. Many thanks also go to the anonymous referee for a number
of helpful comments.
        \end{ack}
\footnotesize


\end{document}